\documentclass[a4paper,11pt]{article}

\usepackage{amsfonts, amsmath, amssymb, amsthm}
\usepackage[colorlinks=true, citecolor=blue]{hyperref}
\usepackage[margin = 2.54cm]{geometry}
\usepackage{mathtools}
\usepackage{enumerate}
\usepackage{verbatim}
\usepackage{tikz}
\usepackage{microtype}
\usepackage{color, tikz, mathtools}
\usepackage{fontenc}
\usepackage[utf8]{inputenc}

\usepackage{tcolorbox}
\usepackage{mdframed}

\usetikzlibrary{positioning}
\usetikzlibrary{shapes}
\usetikzlibrary{decorations.pathreplacing}
\usepackage{enumerate, enumitem}

\newcommand{\NN}{\mathbb{N}}

\newcommand{\RR}{\mathbb{R}}
\newcommand{\R}{\mathbb{R}}
\renewcommand{\SS}{\mathbb{S}}

\newcommand{\cB}{\mathcal{B}}
\newcommand{\cC}{\mathcal{C}}

\newcommand{\cF}{\mathcal{F}}
\newcommand{\cG}{\mathcal{G}}
\newcommand{\cH}{\mathcal{H}}

\newcommand{\cK}{\mathcal{K}}
\newcommand{\cL}{\mathcal{L}}

\newcommand{\cP}{\mathcal{P}}

\newcommand{\cS}{\mathcal{S}}

\newcommand{\ceil}[1]{\left\lceil #1 \right\rceil}

\renewcommand{\epsilon}{\varepsilon}

\newcommand{\dist}{\text{dist}}

\newtheorem{theorem}{Theorem}

\newtheorem{lemma}[theorem]{Lemma}
\newtheorem{claim}[theorem]{Claim}
\newtheorem{definition}[theorem]{Definition}

\newtheorem{observation}[theorem]{Observation}

\newtheorem{proposition}[theorem]{Proposition}

\theoremstyle{definition}

\newcommand{\remove}[1]{}

\begin{document}

\title{
Finite $k$-Transversals of Infinite Families of Fat Convex Sets
}

\author{
	Sutanoya Chakraborty \footnote{Indian Statistical Institute, Kolkata, India}
	\and 
	Arijit Ghosh \footnotemark[1]
	\and 
	Soumi Nandi \footnotemark[1]
}

%\date{}

\maketitle

\begin{abstract}
	We prove an infinite $(p,q)$-theorem for piercing fat compact convex sets in $\RR^d$ with $k$-flats.
	Additionally, we develop a new framework through which infinite $(p,q)$-theorems concerning compact sets and $k$-flats can be extended to their 'colorful' variants.
	Further, we show that the existence of an infinite $(p,q)$-theorem does not necessarily imply a finite $(p,q)$-theorem.
\end{abstract}

\section{Introduction}\label{sec:intro}

Given $d\in\NN$ and a family $\cF$ of sets in $\RR^d$, a collection $\cS$ of sets of $\RR^d$ is a {\em transversal} of $\cF$ if for every $F\in\cF$, there is an $S\in\cS$ such that $F\cap S\neq\emptyset$.
For $k\in\{0,1,\dots,d-1\}$, $\cS$ is called a {\em $k$-transversal} of $\cF$ if $\cS$ is a transversal of $\cF$ and every $S\in\cS$ is a $k$-flat.
When $k=0$ ($k=d-1$), $\cS$ is called a point (hyperplane) transversal of $\cF$. For $p,q\in\NN$ and $p\geq q$, $\cF$ satisfies the {\em $(p,q)$-property} with respect to a family $\cG$ of sets in $\RR^d$ if every subset of $\cF$ of size $p$ contains $q$ sets that are intersected by a set in $\cG$.

Helly's theorem is a well-known result in geometry which states that given a finite family $\cF$ of convex sets in $\mathbb{R}^d$, if the intersection of every $d+1$ sets from the family is nonempty, then there exists a point that lies in every set in the family. In other words, if $\cF$ is a finite family of convex sets satisfying the $(d+1,d+1)$-property with respect to the set of all points in $\RR^d$, then $\cF$ has a point transversal of size $1$.
If $\cF$ contains infinitely many sets, then in addition to the sets being convex, we need the sets to be compact as well for Helly's theorem to hold.
Alon and Kleitman~\cite{AlonK92} proved the following powerful generalization of Helly's theorem:

\begin{theorem}[Alon-Kleitman~\cite{AlonK92}: $(p,q)$-theorem]
    For any three natural numbers $p\geq q\geq d+1$, $\exists c=c(p,q,d)$ such that if $\cF$ is a collection of compact convex sets satisfying the $(p,q)$-property in $\RR^d$ with respect to points, then there exists a finite-sized point transversal of $\cF$.
\end{theorem}
Results of the above kind are called $(p,q)$-theorems, where we have a family $\cF$ of sets in $\RR^d$, another family $\cG$ of sets, and $\cF$ satisfying the $(p,q)$-property with respect to $\cG$ implies the existence of a finite-sized subset $\cS$ of $\cG$ whose size is independent of the size of $\cF$, such that $\cS$ is a transversal of $\cF$.
Alon and Kalai~\cite{AlonK95} proved a $(p,q)$-theorem for hyperplane transversals as well.

For two nonempty sets $A$ and $B$, we say $A$ pierces $B$ if $A\cap B\neq\emptyset$.
Interestingly, there are results demonstrating that the $(p,q)$-theorem fails for $k$-flats piercing compact convex sets of arbitrary shape when $0<k<d-1$.
Alon, Kallai, Mat\u{o}usek and Meshulam \cite{alon2002transversal} showed that there cannot be a $(p,q)$-theorem for piercing compact convex sets with straight lines in $\RR^3$, answering a question Alon and Kalai posed in \cite{AlonK92}.
This indicates that $(p,q)$-theorems involving convex sets and $k$-flats for $k\notin \{0, d-1\}$ would require more conditions than mere compactness.

Along the line of $(p,q)$-properties, Keller and Perles \cite{KellerP22} introduced the notion of $(\aleph_0,q)$-properties, where every infinite subset of some family $\cF$ of sets in $\RR^d$ contains $q\in\NN$ sets that are pierced by an element of some family $\cG$ of sets.
The $(p,q)$-property implies the $(\aleph_0,q)$-property for infinite families, so the latter is less restrictive.
Keller and Perles showed that an $(\aleph_0,k+2)$-theorem can be obtained for $k$-flats piercing $(r,R)$-fat compact convex sets (a compact convex set $C$ is $(r,R)$-fat if $\overline{B}(c,r)\subseteq C\subseteq \overline{B}(c,R)$ for some $c\in C$), $k\in\{0,\dots,d-1\}$, which they later extended to an $(\aleph_0,k+2)$-theorem for $k$-flats piercing compact \textit{near-balls} \cite{KellerP23}.
They defined the notion of \textit{near-balls with parameter $\rho$} as follows: given a $\rho>0$, a family of \textit{near-balls} with parameter $\rho$ is a collection $\cF$ of sets in $\RR^d$ satisfying the following 
property: for every $B\in\cF$, we have $r_B,R_B>0$ satisfying $\frac{r_B}{R_B}\geq\rho$ and $r_B + \rho > R_B$, and there is a $x_B\in B$ such that $\overline{B}(x_B,r_B)\subseteq B\subseteq \overline{B}(x_B,R_B)$.
Jung and P\'alv\"olgyi \cite{jung2024noteinfiniteversionspqtheorems} showed a technique through which $(p,q)$-theorems can be extended to $(\aleph_0,q)$-theorems if the corresponding fractional Helly theorems exist.
They proved a $(p,q)$-theorem for $k$-flats piercing closed balls in $\RR^d$ and fractional Helly theorems for closed balls and $\rho$-fat compact convex sets, and conjectured that a $(p,q)$-theorem holds 
for $k$-flats piercing $\rho$-fat compact convex sets, where given a $\rho>0$, a set $C$ in $\RR^d$ is defined to be $\rho$-fat if there is a 
point $x_C\in C$ and $0<r_C\leq R_C<\infty$, with $\overline{B}(x_C,r_C)\subseteq C\subseteq \overline{B}(x_C,R_C)$, where $\frac{r_C}{R_C}\geq \rho$.\cite{jung2024kdimensionaltransversalsfatconvex}

\section{Our results}\label{sec:results}

Our main result is an $(\aleph_0,k+2)$-theorem for $k$-flats piercing compact $\rho$-fat convex sets.
\begin{theorem}[$(\aleph_0,k+2)$-theorem for $\rho$-fat compact convex families]\label{thm:mono_rho_fat}
    Let $d\in\NN$, $\rho>0$, $\cF$ be an infinite collection of compact convex sets in $\RR^d$ that are $\rho$-fat, and every infinite subset of $\cF$ contains $k+2$ sets that are pierced by a $k$-flat.
    Then there is a finite collection $\cK$ of $k$-flats such that for every $C\in\cF$, there is a $K\in\cK$ with $C\cap K\neq\emptyset$.
\end{theorem}

The arguments we employ in the proof of Theorem \ref{thm:mono_rho_fat} also work for the case where $\cF$ is a family of compact near-balls.
We also demonstrate that the convexity assumption of Theorem \ref{thm:mono_rho_fat} cannot be replaced by a connectedness assumption.
Neither can we do away with the $\rho$-fatness condition on compact convex sets without requiring additional assumptions.
The examples demonstrating the above are in Section \ref{sec:no_go}.\\
We also prove a general framework through which $(\aleph_0,.)$-theorems involving $k$-flats piercing compact sets can be extended to their 'colorful' versions.
To state the result, we need the following definitions of 'colorful':
\begin{definition}[Heterochromatic sequences]\label{def:heterochromatic_seq}
    Let $\{\cF_n\}_{n\in\NN}$ be a sequence of families of sets. 
    A sequence $\{S_n\}_{n\in\NN}$ is a heterochromatic sequence from $\{\cF_n\}_{n\in\NN}$ if there exists a strictly increasing sequence $\{m_n\}_{n\in\NN}$ in $\NN$ such that $S_n\in\cF_{m_n}$.
\end{definition}

The above definition differs from its finite counterpart in that in a heterochromatic sequence chosen from $\{\cF_n\}_{n\in\NN}$, every $\cF_n$ does not need to contribute a member.
Now we define {\em strict heterochromatic sequences}, where every $\cF_n$ must contribute a set.
\begin{definition}[Strict heterochromatic sequences]\label{def:strict_heterochromatic_seq}
    Let $\{\cF_n\}_{n\in\NN}$ be a sequence of families of sets. 
    A sequence $\{S_n\}_{n\in\NN}$ is a strict heterochromatic sequence from $\{\cF_n\}_{n\in\NN}$ if $S_n\in\cF_{n}$.
\end{definition}

\begin{definition}[Heterochromatic and strict heterochromatic $(\aleph_0,q)$-properties]\label{hetero_q_property}
    Let $q\in\NN$ and $\cG$ be a family of sets.
    A sequence $\{\cF_n\}_{n\in\NN}$ of families of sets is said to satisfy the {\em heterochromatic} ({\em strict heterochromatic}) $(\aleph_0,q)$-property with respect to $\cG$ if every {\em heterochromatic} ({\em strict heterochromatic}) sequence from $\{\cF_n\}_{n\in\NN}$ contains $q$ sets that are pierced by an $S\in\cG$.
\end{definition} 
We show the following:
\begin{theorem}[Heterochromatic $(\aleph_0,k+2)$-theorem from monochromatic $(\aleph_0,k+2)$-theorem]\label{thm:mono_to_hetero}
    Let $d\in\NN$, $\cF$ be a collection of compact sets in $\R^d$ and $0\leq k<d$.
    If $\{\cF_n\}_{n\in\NN}$ is a sequence of families of sets where $\cF_n\subseteq \cF$ for all $n\in\NN$ and $\{\cF_n\}_{n\in\NN}$ satisfies the heterochromatic $(\aleph_0,k+2)$-property with respect to $k$-flats, then there is an $m\in\NN$ and a collection $\cK$ of finitely many $k$-flats such that for every $C\in\cF_n$ where $n\geq m$, we have a $K\in\cK$ with $C\cap K\neq\emptyset$.
\end{theorem}

As is evident, the value of $m$ in the above result depends on $\{\cF_n\}_{n\in\NN}$.
This result extends Theorem \ref{thm:mono_rho_fat} to its heterochromatic version.
However, we shall show with an example that replacing {\em heterochromatic sequences} with {\em strict heterochromatic sequences} in Theorem \ref{thm:mono_to_hetero} does not work, even when we take $\cF$ to be the collection of closed unit balls.
\begin{proposition}\label{exmp:balls_hetero_no_go}
    Let $d>1$ be an integer and $k\in\{1,\dots,d-2\}$.
    There exists a sequence $\{\cB_n\}_{n\in\NN}$ of families of closed unit balls in $\RR^d$ with the strict heterochromatic $(\aleph_0,k+2)$-property with respect to $k$-flats, where for every $n\in\NN$, $\cB_n$ does not have a finite $k$-transversal.
\end{proposition}
The strict heterochromatic $(\aleph_0,2)$-property with respect to points also does not guarantee a finite $0$-transversal for any $\cF_n$ when each $\cF_n$ is a collection of closed compact axis-parallel boxes in $\RR^d$, $d>0$. \cite{ChakrabortyGN2023boxesflats}
However, there are still some things that can be said about certain families of sets with the strict heterochromatic $(\aleph_0,2)$-property, as the following result shows.
\begin{theorem}[Strict heterochromatic $(\aleph_0,2)$-theorem for closed unit balls]\label{thm:strong_hetero_balls}
    Let $d\in\NN$ and $\{\cB_n\}_{n\in\NN}$ be families of unit balls in $\RR^d$ satisfying the strict $(\aleph_0,2)$-property with respect to points.
    Then there exist two distinct $n_1,n_2\in\NN$ such that $\cF_{n_1}$ and $\cF_{n_2}$ have a finite point transversal.
\end{theorem}
Clearly, the above result can be extended to families of $(r,R)$-fat convex compact sets.

A natural question in this context is: when can we get an $(\aleph_0,q)$-theorem?
Is the existence of the corresponding $(p,q)$-theorem a necessity?
The answer is no, as the following proposition shows.
\begin{proposition}\label{exmp:near_balls_pq}
    Given $n,q,d\in\NN$, there is a collection $\cF$ of compact near-balls with $|\cF|=n$, satisfying the $(q,q)$-property with respect to points, that cannot be pierced by $\Theta(n)$ points.
\end{proposition}

\section{Organization of the paper}

In Section \ref{sec:rho_fat}, we prove Theorem \ref{thm:mono_rho_fat}. The proofs of Theorems \ref{thm:mono_to_hetero} and \ref{thm:strong_hetero_balls} are in Section \ref{sec:heterochromatic} (Subsections \ref{ssec:mono_to_hetero} and \ref{ssec:strict_hetero}, respectively), and Propositions \ref{exmp:balls_hetero_no_go} and \ref{exmp:near_balls_pq} are proved in Section \ref{sec:no_go} (Subsections \ref{ssec:limits_strict_hetero} and \ref{ssec:pq_no_go}, respectively).

\section{Definitions and notations}
This section contains definitions and notations that are not described in Sections \ref{sec:intro} or \ref{sec:results}.
\begin{itemize}
    \item 
        For any $n \in \NN$, $[n]$ denotes the set $\left\{ 1, \dots, n\right\}$.

    \item 
        The origin in $\mathbb{R}^{d}$ is denoted by $O$. 

    \item 
        $\mathbb{S}^{d-1}$ denotes the unit sphere in $\RR^d$ with center at the origin $O$.

    \item 
        For any $A\subseteq\RR^d$, $\overline{A}$ denotes the closure of $A$ in the standard topology of $\RR^d$.

    \item 
        For any $p\in\RR^d$ and $\epsilon>0$, $B(p,\epsilon)$ denotes the open ball of radius $\epsilon$ centered at $p$, and $\overline{B}(p,\epsilon)$ denotes the closure of $B(p,\epsilon)$.
        
    \item 
        For any $x\in\RR^d$, $||x||$ denotes the $L_2$-norm of $x$.

    \item 
        For any $x\in\RR^d$ and a set $S\subseteq\RR^d$, $\dist(x,S) : = \inf\{\dist(x,b)\;|\;b\in B\}$. For two sets $A$ and $B$ in $\RR^d$, $\dist(A,B)=\inf\{||a-b||\mid a\in A,\; b\in B\}$.
    
    \item 
        For $k\in\{0\leq k\leq d-1\}$, a $k$-dimensional affine subspace of $\RR^d$ is called a \emph{$k$-flat}.
        A point is a $0$-flat, a line is a $1$-flat, and a hyperplane in $\RR^d$ is $(d-1)$-flat.

    \item 
        Given a $p\in\RR^d$, a ray $\vec{r}$ from $p$ is a set $\{p+tv\mid t\geq 0\}$ where $v\in\SS^{d-1}$.

      \item Given any set $U\in\mathbb{R}^d$, $\text{pos}(U):=\{cu\mid c\geq 0, u\in U\}$.

    \item 
      Given any $\epsilon>0$ and a ray $\vec{r}=\{p+tv\mid t\geq 0\}$ from $p$ in $\RR^d$, $v\in\SS^{d-1}$, $\text{cone}(\vec{r},\epsilon)=\{p+u\mid u\in \text{pos}(\overline{B}(v,\epsilon)\}$. 

    \item 
      Given a ray $\vec{r}=\{p+tv\mid t\geq 0\}$, $v\in \SS^{d-1}$ and $n\in\NN$, we write $R(\vec{r},\frac{1}{n}) := \{p+u\mid u\in \text{pos}(\overline{B}(v,\frac{1}{n}))\cap \overline{B}(p,\frac{1}{n})$ and 
    $Q(\vec{r},n) := \{p+u\mid u\in\text{pos}(B(v,\frac{1}{n}))\}\setminus B(p,n)$.

    \item 
        For a family $\cF$ of sets in $\RR^d$ and $k\in\{0,\dots,d-1\}$, a collection $\cK$ of $k$-flats is a {\em $k$-transversal} of $\cF$ if for every $S\in\cF$, there is a $K\in\cK$ with $S\cap K\neq\emptyset$.
        $\cF$ is said to have a finite $k$-transversal if it has a $k$-transversal of finite size.

    \item 
        Given $k\in\{0,\dots,d-1\}$, a family $\cS$ of sets in $\RR^d$ is $k$-independent if there is no $k$-flat piercing $k+2$ distinct sets in $\cS$.

    \item\emph{$k$-set:}
    A family $\mathcal{F}$ of sets in $\RR^d$ is a $k$-set if $\cF$ does not have a finite $k$-transversal.

    \item For a family $\cF$ of sets in $\R^d$ and any two positive numbers $r,R\in\R^{+}$, $\cF_{r,R}:=\{S\in\cF\mid r_S\geq r \text{ and }R_S\leq R\}$, $\cF_{\leq r}:=\{S\in\cF\mid r_S\leq r\}$ and $\cF_{\geq r}:=\{S\in\cF\mid r_S\geq r\}$.

    \item Given $\rho>0$, for a family $\cF$ of $\rho$-fat sets, let for every $S\in\cF$, $x_S\in S$ denote a point for which there exist $r_S,R_S>0$ such 
      that $\overline{B}(x_S,r_S)\subseteq S\subseteq \overline{B}(x_S,R_S)$, $\frac{r_S}{R_S}\geq\rho$.
      If $S$ is also a near-ball with parameter $\rho$, then $x_S$ denotes a point in $S$ such that $\overline{B}(x_S,r_S)\subseteq S\subseteq \overline{B}(x_S,R_S)$, $\frac{r_S}{R_S}\geq\rho$,
      $\rho+r_S\geq R_S$.

\end{itemize}

\section{An $(\aleph_0,k+2)$-theorem for $\rho$-fat compact convex sets}\label{sec:rho_fat}

In this section, we prove Theorem \ref{thm:mono_rho_fat}.
The arguments employed in the proof readily extend to the case involving compact near-balls as well, so we shall prove the following theorem from which Theorem \ref{thm:mono_rho_fat} follows immediately:

\begin{theorem}\label{thm:rho_fat_near_ball}
    Let $d\in\NN$, $\rho>0$ and $\cF$ be a family of compact sets that are either convex and $\rho$-fat or are near-balls with parameter $\rho$.
    Then either $\cF$ has a finite $k$-transversal, or $\cF$ contains an infinite $k$-independent subset.
\end{theorem}

We make the following trivial observation below:

\begin{observation}
    For any $k$-set $\cF$, at least one of the following must be true:
    \begin{itemize}
        \item [(1)] For every $r>0$, $\cF_{\leq r}:=\{S\in\cF\;|\; r_S<r\}$ is a $k$-set.

         \item [(2)] There exists an $r>0$ for which $\cF_{\geq r}:=\{S\in\cF\;|\; r_S\geq r\}$ is a $k$-set.
    \end{itemize}
\label{obs2}\end{observation}

\begin{observation}\label{obs: bddregion}
    Let $\cF$ be a $k$-set of either either compact convex $\rho$-fat sets or compact near-balls with parameter $\rho>0$, such that there is a bounded set $A\subset \RR^d$ and $\overline{A}$ intersects every set in $\cF$.
    Then for every $r>0$, $\cF_{\geq r}$ is not a $0$-set.
\end{observation}

\begin{proof}
    Here we show that for every $r>0$, $\cF_{\geq r}$ is pierceable by finitely many points.
    Recall from the definitions that for every $S\in\cF_{\geq r}$, $B_S$ denotes a largest inscribed ball in $S$.
    Now, let $S$ be a compact near-ball in $\cF$ that pierces $\overline{A}$.
    Then, if $\overline{A}\subset B(O,n)$ for some $n\in\mathbb{N}$, then $B_S\cap B(O,n+\rho)\neq\emptyset$.
    It is easy to see that a finite number of points suffice to pierce all closed balls of radius $\geq r$ that pierce $B(O,n+\rho)$.\\
    It is easy to see that there are a finite number of points $P\subset\RR^d$ such that for every $S\in\cF_{\geq r}$ with $B_S\cap \overline{A}\neq\emptyset$, we have $S\cap P\neq\emptyset$.
    That is because the radius of $B_S$ is at least $r$, and all these $B_S$ intersect $\overline{A}$.
    Now, let $\cS\subseteq\cF_{\geq r}$ be the collection of all compact convex sets $S$ in $\cF_{\geq r}$ for which $B_S$ does not intersect $\overline{A}$.
    Then there is an $\epsilon>0$ such that for every $S\in\cS$, there is a point $p_S$ on $\partial A\cap S$ from which a finite cone $C_S$ of aperture $\epsilon$ and axis length $r$ exists that lies completely within $S$.
    It is easily seen that we can find a finite set $Q\subset \R^d$ such that for all $S\in\cS$, $Q\cap C_S\neq \emptyset$, since for any two $S_1,S_2\in\cS$, $C_{S_1}$ and $C_{S_2}$ are identical up to translation and rotation, every $C_S$ for $S\in\cS$ contains a ball with volume greater than $\beta r^d$ for some $\beta>0$, and there is a bounded region that contains $C_S$ for every $S\in\cS$. 
\end{proof}

\begin{claim}\label{claim1}
        Given a $k$-set $\cF$ of either compact convex $\rho$-fat sets, or compact near-balls with parameter $\rho>0$, there is a point $p\in\RR^d$ and a ray $\vec{r}$ from $p$ such that one of the following two cases must arise:
    \begin{itemize}
      \item\label{ite:point_cluster_rhofat} [(1)] $\forall n\in\NN,\;\cF_n:=\{S\in\cF\;|\;x_S\in R(\vec{r},\frac{1}{n})\}$ is a $k$-set.

      \item\label{ite:point_faraway_rhofat} [(2)] $\forall n\in\NN,\;\cF_n:=\{S\in\cF\;|\;x_S\in Q(\vec{r},n)\}$ is a $k$-set.
    \end{itemize}
\end{claim}

\begin{proof}
    Suppose $\cG_n:=\{S\in\cF\;|\;x_S\in\RR^d\setminus B(O,n)\}$. Then either $\forall n\in\NN,\; \cG_n$ is a $k$-set or $\exists n'\in\NN$ such that $\cG_{n'}$ is not a $k$-set.

\vspace{10pt}

\noindent
{\bf Case I: $\forall n\in\NN,\; \cG_n$ is a $k$-set.} First choose $\epsilon>0$ such that $\epsilon<<1$. Now we consider a finite cover of $\SS^{d-1}$ with open balls of radius $\epsilon$ and center on $\SS^{d-1}$. We call an open set $U$ in $\RR^d$ \emph{good} with respect to $\cG_n$ if $\cG_{B,n}:=\{S\in\cG_n\;|\;\frac{x_S}{||x_S||}\in U\cap\SS^{d-1}\}$ is a $k$-set. Let $B_{1,1},\dots, B_{1,t_1}$ be the good sets with respect to $\cG_1$ chosen from the finite cover of $\SS^{d-1}$ of open balls with radius $\epsilon$ and center on $\SS^{d-1}$. Let 
\begin{align*}
        V_{1} = \left( \bigcup_{i=1}^{t_1} \Bar{B}_{1,i} \right) \bigcap \SS^{d-1}.
\end{align*}
    For $\cG_n$, we set 
    $$
        V_{n} = \left( \bigcup_{i=1}^{t_n} \Bar{B}_{n,i} \right) \bigcap \SS^{d-1},
    $$
    where $B_{n,i}$s are the good sets with respect to $\cG_n$ chosen from a finite cover of $\SS^{d-1}$ of open balls with radius ${\epsilon}/{n}$ and center on $\SS^{d-1}$.

For each $n\in\NN$, we have 
$\bigcap_{i=1}^n V_i\neq\emptyset$.
To see why, note that if an open set $U\subset\RR^d$ is good with respect to $\cG_{n}$, then for any $i\in[n]$, in a finite cover of $U\cap \SS^{d-1}$ of open balls of radius $\frac{\epsilon}{i}$, there is an open ball $B_{i,j}$ such that $U_{i}=B_{i,j}\cap U$ is good with respect to $\cG_{n}$.

Since each $V_i$ is compact, this implies that 
$\bigcap_{n\in\NN}V_n\neq\emptyset$. 
Now choose $v\in\bigcap_{n\in\NN}V_n$. Let $\vec{r}:=\{cv\;|\;c\geq0\}$. Then we get that, for each $n\in\NN$, $\cF_{n} = \left\{S\in\cF\;|\;x_S\in Q(\vec{r},n)\right\}$ is a $k$-set.

\vspace{10pt}

\noindent
{\bf Case II: $\exists n'\in\NN$ such that $\cG_{n'}$ is not a $k$-set.}
We claim that $\exists p\in\Bar{B}(O,n')$ such that for every open ball $B$ with $p\in B$, we have $\cF_B=\{S\in\cF\;|\;x_S\in B\}$ to be a $k$-set. Otherwise, for every $a\in \Bar{B}(O,n')$, we have an open ball $B_a$ with  $a\in B_a$ such that $\cF_{B_a}$ is not a $k$-set. 
By compactness of $\Bar{B}(O,n')$, we get that the open cover $\{B_a\;|\;a\in\Bar{B}(O,n')\}$ of $\Bar{B}(O,n')$ has a finite subcover $B_{a_1},\dots,B_{a_t}$. 
This implies that $\cG'_{n'}=\{S\in\cF\;|\;x_S\in \overline{B}(O,n')\subset \cup_{i=1}^{t} B_{a_i}\}$ is not a $k$-set. 
By our assumption, $\cG_{n'}$ is not a $k$-set, and therefore, $\cF=\cG_{n'}\cup\cG'_{n'}$ is not a $k$-set. 
But this contradicts the fact that $\cF$ is a $k$-set. 

Now we can safely assume that for every $S\in\cF,\; x_S\neq p$, and without loss of generality, we assume that $p$ is the origin $O$. 
Next, as in Case $I$, pick $\epsilon>0$ with $\epsilon<<1$ and cover $\SS^{d-1}$ with finitely many open balls of radius $\epsilon$ and center on $\SS^{d-1}$. 
For each $n\in\NN$, let $\cH_{n} := \left\{S\in\cF\;|\;x_S\in B(O,\frac{1}{n}) \right\}$. 
Without loss of generality, we assume that $\forall S\in\cF, x_S\neq O$. 
We call an open set \emph{good} with respect to $\cH_n$ if 
$$
    \cH_{B,n} := \left\{S\in\cF\;|\;\frac{x_S}{||x_S||}\in B\cap\SS^{d-1}\right\}
$$ 
is a $k$-set. Let $B_{1,1},\dots, B_{1,t_1}$ be the good sets with respect to $\cH_1$ chosen from a finite cover of $\SS^{d-1}$ of open balls with radius $\epsilon$ and center on $\SS^{d-1}$. 
Define $V_{1} =\left( \bigcup_{i=1}^{t_1}\Bar{B}_{1,i} \right)\cap\SS^{d-1}$. Similarly, for each $n\in\NN$, we define $V_{n} = \left( \bigcup_{i=1}^{t_n}\Bar{B}_{n,i} \right)\cap\SS^{d-1}$, where $B_{n,1},\dots, B_{n,t_n}$ are the good sets with respect to $\cH_n$ chosen from a finite open cover of $\SS^{d-1}$ of open balls with radius $\frac{\epsilon}{n}$ and center on $\SS^{d-1}$. 
As argued in Case-I, $\bigcap_{n\in\NN}V_n\neq\emptyset$. Let $v\in\bigcap_{n\in\NN}V_{n}$ and $\vec{r}:=\{cv\;|\;c\geq 0\}$. 
Then we get that for each $n\in\NN$, $\cF_n:=\{S\in\cF\;|\;x_S\in R(\vec{r},\frac{1}{n})\}$ is a $k$-set. 
\end{proof}

Now we prove Theorem \ref{thm:rho_fat_near_ball} when $k=0$.

\begin{theorem}[$(\aleph_0,2)$-theorem for compact convex $\rho$-fat sets and compact near-balls with parameter $\rho$ and points]
  Let $\rho\geq 1$ and $\mathcal{F}$ be a family of compact $\rho$-fat sets such that if $S\in\mathcal{F}$, then either $S$ is convex or 
  a near-ball with parameter $\rho$.
  Then if $\mathcal{F}$ is a $0$-set, then there is an infinite sequence $\{S_n\}_{n\in\mathbb{N}}\subset\mathcal{F}$
  that is $0$-independent, i.e., $S_i\cap S_j=\emptyset$ whenever $i\neq j$, $i,j\in\mathbb{N}$.
  \label{thm:rho_fat_0}
\end{theorem}

\begin{proof}
  If for every bounded set $U$ in $\mathbb{R}^d$, the collection $\{S\in\mathcal{F}\mid 
  S\cap U=\emptyset\}$ is an $0$-set, then we can easily find a sequence $\{S_n\}_{n\in\mathbb{N}}$
  from $\mathcal{F}$ that is $0$-independent, where after choosing $S_1,\dots,S_m$,
  we can pick $S_{m+1}\in\mathcal{F}$ such that it does not intersect $\cup_{j\in[m]}S_j$, as 
  $\cup_{j\in[m]}S_j$ is a bounded set.
  So let there be a bounded set $U$ such that $\{S\in\mathcal{F}\mid S\cap U=\emptyset\}$
  is not a $0$-set.
  Then, without loss of generality, we can assume that every $S\in\mathcal{F}$ intersects $U$.

  \noindent\textbf{Case-I:} Claim \ref{claim1}, (1) holds for some $p\in\mathbb{R}^d$ and $\vec{r}$.
  \newline Let $\forall n\in\mathbb{N}$, $\mathcal{F}_n:=\{S\in\mathcal{F}\mid x_S\in R(\vec{r},\frac{1}{n})\}$.
  Observation \ref{obs: bddregion} shows that for every $r>0$ $\mathcal{F}_{\leq r}\cap\mathcal{F}_n$ is 
  a $0$-set.
  Choose $S_1\in\mathcal{F}_1$ such that $S_1$ does not contain $p$.
  Then there is a closed ball $\overline{B}(p,r_2)$ such that $S_1\cap \overline{B}(p,r_2)=\emptyset$, $r_2>0$.
  We have an $n_2\in\mathbb{N}$ such that $r_2>(1+2\rho)\frac{1}{n_2}$.
  As $\mathcal{F}_{\leq \frac{1}{n_2}}\cap \mathcal{F}_{n_2}$ is a $0$-set, we can choose a 
  $S_2\in\mathcal{F}_{\leq \frac{1}{n_2}}\cap \mathcal{F}_{n_2}$ such that $S_2$ does not 
  contain $p$. 
  This implies that for any $y\in S_2$, $\dist(y,p)\leq \dist(y,x_S)+\dist(x_S,p)< (1+2\rho)\frac{1}{n_2}<r_2$.
  Therefore, $S_2\cap S_1=\emptyset$.
  In general, given a finite sequence $S_1,\dots,S_m$ of mutually disjoints sets from $\mathcal{F}$
  such that $p\notin S_j$ for every $j\in[m]$, we can find an $r_{m+1}>0$ such that $\overline{B}(p,r_{m+1})\cap S_j=\emptyset$
  for every $j\in[m]$, and proceed as before to obtain $S_{m+1}$ from $\mathcal{F}$.
  Continuing this way, we obtain an infinite sequence $\{S_n\}_{n\in\mathbb{N}}$ in $\mathcal{F}$ with $S_i\cap S_j=\emptyset$ 
  whenever $i\neq j$, $i,j\in\mathbb{N}$.

  \noindent\textbf{Case-II:} Claim \ref{claim1}, (2) holds for some $\vec{r}$. 
  \newline Then Observation \ref{obs: bddregion} tells us that for every $r>0$, 
  $\mathcal{F}_{\leq r}\cap\mathcal{F}_n$ is an $0$-set, where $\mathcal{F}_n:=\{S\in\mathcal{F}\mid x_S\in Q(\vec{r},n)\}$.
  Without loss of generality, let $S\cap\vec{r}=\emptyset$ for every $S\in\mathcal{F}$, as $\mathcal{F}$ 
  is a $0$-set.
  Choose any $S_1\in\mathcal{F}_1$.
  Then there is an $n_2\in\mathbb{N}$ such that $\dist(S_1,Q(\vec{r},n_2))=r'_2>0$.
  Let $r_2>0$ be a positive number satisfying $(2\rho+1)r_2 < r'_2$.
  Then if $S_2\in\mathcal{F}_{\leq r_2}\cap\mathcal{F}_{n_2}$, then for any point 
  $y\in S_2$, $\dist(S_1,S_2)\geq \dist(x_S,S_1)-\dist(x_S,y)>r'_2-2\rho r_2>0$.
  Therefore, $S_1\cap S_2=\emptyset$.
  In general, given a finite sequence $S_1,\dots,S_m$ of mutually disjoints sets from $\mathcal{F}$,
  we can find an $n_{m+1}\in\mathbb{N}$ for which $Q(\vec{r},n_{m+1})\cap S_i=\emptyset$ 
  for every $i\in[m]$.
  Then there is an $r'_{m+1}>0$ such that $\dist(S_i,Q(\vec{r},n_{m+1}))>r'_{m+1}$ for 
  every $i\in[m]$, and we proceed as before to obtain $S_{m+1}$ from $\mathcal{F}$.
  Continuing this way, we obtain an infinite sequence $\{S_n\}_{n\in\mathbb{N}}$ in $\mathcal{F}$ with $S_i\cap S_j=\emptyset$ 
  whenever $i\neq j$, $i,j\in\mathbb{N}$.
\end{proof}

\begin{claim}
    Let $\rho>0$ and $\cF$ be a $k$-set of either compact convex $\rho$-fat sets, or compact near-balls with parameter $\rho$, where $0<k\leq d-1$.
    Suppose $p\in\RR^d$ and $\vec{r}$ is a ray from $p$ such that $\forall n\in\NN$, $\{\cF_n\}_{n\in\NN}$ is a nested sequence of $k$-sets with $\cF_n\subseteq\cF$, and for all $S\in\cF_n$, we have $x_S\in R(\vec{r},\frac{1}{n})$.
    Then $\cG_n:=\{S\in\cF_n\;|\;S\subset R(\vec{r},\frac{1}{n})\}$ is a $k$-set. 
\label{claim2}\end{claim}

\begin{proof}
   Without loss of generality, we assume that $p$ is the origin $O$ and $\vec{r}$ is the positive $d$-th axis, i.e, $\vec{r}=\{ce_d\in\RR^d\;|\;c\geq 0\}$, where $e_d=(0,\dots,0,1)$. 
   We also assume that $\forall S\in\cF,\;\vec{r}\cap S=\emptyset$, as $k>0$.
   Now we define a function $m:\NN\rightarrow\NN$ by $m(n)=\ceil{2\rho(2^n+2\rho)}$ and pick any $n\in\NN$. 
   We choose any $S\in\cF_{m(n)}$. 
   Then $r_S<\dist(x_S,\vec{r})$. 
   Now we take the $u\in\vec{r}$ for which $\dist(x_S,\vec{r})=\dist(x_S,u)$. Clearly, $\frac{||x_S-u||}{||u||}\leq\frac{1}{m(n)}$. Then 
   $$
        \frac{r_S}{||u||}\leq\frac{1}{m(n)}\implies \frac{2R_S}{||u||}\leq\frac{2\rho}{m(n)}\leq\frac{1}{2^n}.
   $$ 
   Also,
   $$
        \frac{2R_S}{||u||-2R_S}\leq\frac{2\rho}{m(n)-2\rho}\leq\frac{1}{2^n}.
   $$ 
   Pick any $y=(y_1,\dots,y_d)\in S$.
   Then $y_d>||u||-2R_S$, because $S$ is fully contained in $B(x_S,2R_S)$, and the $d$-th coordinate of $x_S$ is $||u||$.
   Also, the closest point on $\vec{r}$ to $y$ is $y_d e_d$, and $\dist(y,\vec{r})=\sqrt{y_1^2+\dots+y_{d-1}^2}<2R_S$.
   So, $||\frac{y}{y_d}-e_d||<\frac{2R_S}{||u||-2R}\leq \frac{1}{2^{n}}$.
   Therefore, $\frac{y}{y_d}\in B(e_d,\frac{1}{2^n})\subset B(e_d,\frac{1}{n})$.
   Also, $||y||<4R_S + ||u||\leq ||u||(\frac{1}{2^{n-1}}+1)\leq\frac{1}{m(n)}(\frac{1}{2^{n-1}}+1)<\frac{1}{n}$.
   This shows that $S\subset R(\vec{r},\frac{1}{n})$. So $\cG_n\supseteq\cF_{m(n)}$ and hence $\cG_n$ is a $k$-set.   
\end{proof}

Now we prove Theorem \ref{thm:rho_fat_near_ball} for the case where $\forall n\in\NN,\;\cF_{n} := \left\{S\in\cF\;|\;x_S\in R(\vec{r},\frac{1}{n})\right\}$ is a $k$-set.

\begin{lemma}\label{lem1}
       Let $\rho>0$ and $\cF$ be a family of either compact convex $\rho$-fat sets or compact near-balls with parameter $\rho$ in $\RR^d$. 
       If there is a point $p\in\RR^d$ and a ray $\vec{r}$ from $p$ such that $\forall n\in\NN,\;\cF_n:=\{S\in\cF\;|\;x_S\in R(\vec{r},\frac{1}{n})\}$ is a $k$-set, then $\cF$ contains a $k$-independent sequence.
\end{lemma}

\color{black}\begin{proof} 
We shall prove the result by induction on $k$.

    \paragraph{Induction hypothesis:} 
    Given any $d\in\mathbb{N}$, let $p\in\mathbb{R}^d$ be a point such that $\forall n\in\NN, \cF_n$ 
    is a $(k-1)$-set of either convex compact $\rho$-fat sets or compact near-balls with parameter $\rho$, 
    $\cF_n\supseteq\cF_{n+1}$, and $\forall S\in\cF_n$, we have $x_S\subset B(p,\frac{1}{n})\}$.
    Then there is a $(k-1)$-independent sequence $\{S_n\}_{n\in\NN}$ with $S_n\in\cF_n$ for each 
    $n\in\NN$ satisfying $p\notin\cC_k(S_{i_1},\dots,S_{i_k})$, where $\cC_k(S_{i_1},\dots,S_{i_k})$ is the union of all $(k-1)$-flats passing through $S_{i_1},\dots,S_{i_k}$.

     \paragraph{Inductive step:}
    The base case is proven in Theorem \ref{thm:rho_fat_0}.
    So let $k>0$.
    By Claim~\ref{claim1} and Claim~\ref{claim2}, there exists a $p\in\mathbb{R}^d$ and a ray $\vec{r}$ such that $\forall n\in\NN,\;\cG_n:=\{S\in\cF_n\;|\;S\subset R(\vec{r},\frac{1}{n})\}$ is a $k$-set. 
    Without loss of generality, we assume that $p$ is the origin $O$ and $\vec{r}$ is the positive $d$-th co-ordinate axis. 
    Let $Y$ be the hyperplane perpendicular to $\vec{r}$ passing through $O$. 
    Then by the induction hypothesis there exists a sequence $\{S_n\}_{n\in\NN}$ with $S_n\in\cF_n$ 
    such that the sequence of the projections of each $S_n$ on $Y$ is a $(k-1)$-independent sequence and no $(k-1)$-flat 
    on $Y$ passing through the projections of $k$ distinct $S_n$'s contains the origin $O$.   
    Further, we have that for every $n\in\NN$, $S_n$ is contained in $R(\vec{r},\frac{1}{n})$.\\ 
    We proceed to select our sequence in the following manner: we pick $k$ distinct sets $S_{i_1},\dots, S_{i_k}$ from $\{S_n\}_{n\in\NN}$. Consider the union of all $k$-flats passing through the origin $O$ and $S_{i_1},\dots, S_{i_k}$, denoted by $\cC_{k+1}(S_{i_1},\dots, S_{i_k},\{O\})$. Then we can show that $\cC_{k+1}(S_{i_1},\dots, S_{i_k},\{O\})\cap\vec{r}=\{O\}$. That is because if there is a $k$-flat through $S_{i,1},\dots,S_{i_k}$ that contains $\vec{r}$, then its projection onto $Y$ is a $(k-1)$-flat that passes through $O$, which is a contradiction.
    This means $\exists n_{k+1}\in\NN$ such that $\cC_{k+1}(S_{i_1},\dots, S_{i_k},\{O\})\cap R(\vec{r},\frac{1}{n_{k+1}})=\{O\}$. We pick $S_{i_{k+1}}\in\{S_n\}_{n\in\NN}$ such that $i_{k+1}>i_j\;\forall j\in [k]$ and $S_{i_{k+1}}\subset R(\vec{r},\frac{1}{n_{k+1}})$. So $O\not\in \cC_{k+1}(S_{i_1},\dots, S_{i_k},S_{i_{k+1}})$.
    Define $\cS_{k+1}=\{S_{i_1},\dots,S_{i_{k+1}}\}$.

    Now let $\cS_m=\{S_{i_1},\dots, S_{i_m}\}\subset \{S_n\}_{n\in\NN}$ with $i_1<i_2<\dots <i_m$ such that for any distinct $k+1$ sets $S'_1,\dots, S'_{k+1}$ chosen from $S_{i_1},\dots,S_{i_m}$ we have $O\not\in\cC_{k+1}(S'_1,\dots,S'_{k+1})$ and for any $k+2$ distinct sets $S'_1,\dots, S'_{k+2}$ chosen from $S_{i_1},\dots,S_{i_m}$ we have $\cC_{k+1}(S'_1,\dots,S'_{k+1})\cap  S'_{k+2}=\emptyset$.
    We will show that we can extend $\cS_m$.
    Since for any distinct $k+1$ sets $S'_1,\dots, S'_{k+1}\in \cS_m$, $O\not\in\cC_{k+1}(S'_1,\dots,S'_{k+1})$, there exists $n_1\in\NN$ such that $B(O,n_1)\cap\cC_{k+1}(S'_1,\dots,S'_{k+1})=\emptyset$ for every distinct $S'_{1},\dots,S'_{k+1}\in\cS_m$.
    Also for any distinct $k$ sets $S'_1,\dots, S'_{k}$ chosen from $\cS_m$ we have $\exists n_2\in\NN$ such that $R(\vec{r},\frac{1}{n_2})\cap\cC_{k+1}(S'_1,\dots,S'_{k},\{O\})=\{O\}$. Now we pick $S_{i_{m+1}}\in\{S_n\}_{n\in\NN}$ such that $i_{m+1}>i_{m}$ and $S_{i_{m+1}}\subset R(\vec{r},\frac{1}{\Bar{n}})\setminus\{O\}$, where $\Bar{n}=\max\{n_1,n_2\}$. 
    Now, define $\cS_{m+1}=\cS_m\cup\{S_{i_{m+1}}\}$, which is our required extension of $\cS_m$.
    Continuing this way, we obtain our required $k$-independent sequence.
\end{proof}

Next we prove the case when $\forall n\in\NN,\;\cF_n:=\{S\in\cF\;|\;x_S\in Q(\vec{r},n)\}$ is a $k$-set. 
Before that, we require the following result, similar in spirit to Claim \ref{claim2}.
\begin{claim}
    Let $\rho>0$ and $\cF$ be a $k$-set in $\RR^d$ of either compact convex $\rho$-fat sets or compact near-balls with parameter $\rho$, where $0<k\leq d-1$.
    Suppose $p\in\RR^d$ and $\vec{r}$ is a ray from $p$ such that $\forall n\in\NN$, $\{\cF_n\}_{n\in\NN}$ is a nested sequence of $k$-sets with $\cF_n\subseteq\cF$, and for all $S\in\cF_n$, we have $x_S\in Q(\vec{r},n)$.
    Then $\forall n\in\NN,\;\cG_n=\{S\in\cF_n\;|\:S\subset Q(\vec{r},n)\}$ is a $k$-set.
    \label{claim3}
\end{claim}

\begin{proof}
Without loss of generality, let $p$ be the origin $O$ and $\vec{r}$ be the positive $d$-th axis. 
Also assume that $\forall S\in\cF,\; S\cap\vec{r}=\emptyset$, since $k>0$.
Define, $\forall n\in\NN,\;m(n)=\ceil{2\rho(2^n+2\rho)}$.

Choose any $S\in\cF_{m(n)}$. 
Let $u\in\vec{r}$ be such that $\dist(x_S,\vec{r})=\dist(x_S,u)$. Since $S\cap\vec{r}=\emptyset$, we have $r_S\leq ||x_S-u||$. 
Also 
$$
    \frac{||x_S-u||}{||u||}\leq\frac{1}{m(n)}
    \implies\frac{r_S}{||u||}\leq\frac{1}{m(n)}\implies\frac{2R_S}{||u||}\leq\frac{2\rho}{m(n)}\leq\frac{1}{2^n}.
$$
Again, observe that 
$$
    \frac{2R_S}{||u||-2R_S}\leq\frac{1}{2^n}.
$$ 
Let $y=(y_1,\dots,y_d)\in S$.
Arguments similar to those in Claim \ref{claim2} show that $\frac{y}{y_d}\in B(e_d,\frac{1}{n})$.
Also, $||y||>||u||-2R_S\geq ||u||(1-\frac{1}{2^n})\geq m(n)(1-\frac{1}{2^n})>n$.
This shows that $\cF_{m(n)}\subseteq \cG_n$ and hence $\cG_n$ is a $k$-set.
\end{proof}
Now we prove the following lemma, which will conclude the proof of Theorem~\ref{thm:rho_fat_near_ball}.

\begin{lemma}\label{lemQ}
       Let $\rho>0$ and $\cF$ be a family of either compact convex $\rho$-fat sets or compact near-balls with parameter $\rho$.
       If there is a point $p\in\RR^d$ and a ray $\vec{r}$ from $p$ such that $\forall n\in\NN,\;\cF_n:=\{S\in\cF\;|\;x_S\in Q(\vec{r},n)\}$ is a $k$-set, then $\cF$ contains a $k$-independent sequence.
\end{lemma}

\begin{proof}
    We induct on $k$.

    \paragraph{Induction hypothesis:} 
    Given any $d\in\mathbb{N}$, let, $\forall n\in\NN,\;\cF_n$ be a $(k-1)$-set of either compact convex $\rho$-fat sets or compact near-balls with parameter $\rho$
    such that $\cF_n\supseteq\cF_{n+1}$ and $\forall S\in\cF_n$, we have $x_S\subset \mathbb{R}^{d}\setminus B(O,n)$. 
    Then there is a $(k-1)$-independent sequence $\{S_n\}_{n\in\mathbb{N}}$ such that $S_n\in\cF_n$. 

    \paragraph{Inductive step:} 
    The base case is proven in Theorem \ref{thm:rho_fat_0}.
    Now we induct on $d$ and $k$. 
    We have already proved the case for arbitrary $d\in\NN$ and $k=0$, so we assume that $k>0$.
    By Lemma~\ref{lem:finite_piercing}, we can assume that $\cF$ is a countably infinite set.
    By Claim~\ref{claim3}, there exists a $p\in\mathbb{R}^d$ and a ray $\vec{r}$ from $p$ such that $\cG_n=\{S\in\cF_n\mid S\subset Q(\vec{r},n)\}$ is a $k$-set.
    Without loss of generality, we can assume that $\forall n\in\mathbb{N}$, $\cF\setminus\cG_{n}$ is a finite set, since for each $n\in\mathbb{N}$, we can pick a finite subset from $\cG_n$ that cannot be pierced by $n$ $k$-flats and consider $\cF$ to be the union of all these subsets.
    Further, let $p$ be the origin $O$, and let $Y$ be the hypeprlane through $O$ perpendicular to $\vec{r}$.
    Let $\cF'$ consist of the projections of all $S\in\cF$ onto $Y$.
    Then $\cF'$ is a $(k-1)$-set, and therefore, Claim~\ref{claim1} applies on $\cF'$.
    If Claim~\ref{claim1}~(1) is satisfied, then Lemma~\ref{lem1} applies, and if Claim~\ref{claim1}~(2) is satisfied, then our induction hypothesis applies, so in either case 
    we get a sequence $\{S_n\}_{n\in\mathbb{N}}$ in $\cF$ whose projection on $Y$ is $(k-1)$-independent.

    Now pick any distinct $k+1$ sets $S_{i_1}, S_{i_2},\dots, S_{i_{k+1}}$ from $\{S_n\}_{n\in\NN}$. Let $\cC_{k+1}(S_{i_1}, \dots, S_{i_{k+1}})$ denote the union of all $k$-flats passing through $S_{i_1}, S_{i_2},\dots, S_{i_{k+1}}$. 
    Clearly, there is no $k$-flat that passes through $S_{i_1}, S_{i_2},\dots, S_{i_{k+1}}$ that contains a straight line parallel to $\vec{r}$.
    We claim that if every $k$-flat passing through $S_{i_1}, S_{i_2},\dots, S_{i_{k+1}}$ does not contain any straight line parallel to $\vec{r}$, then $\exists n'\in\NN$ such that $Q(\vec{r},n')\cap\cC_{k+1}(S_{i_1}, S_{i_2},\dots, S_{i_{k+1}})=\emptyset$. We shall prove this as Claim~\ref{claim4}.
    Since there is a $S_{i_{k+2}}\subset Q(\vec{r},n')$, we have $S_{i_{k+2}}\cap\cC_{k+1}(S_{i_1}, S_{i_2},\dots, S_{i_{k+1}})=\emptyset$. 
    In general, for $S_{i_1},\dots,S_{i_m}$ chosen from $\{S_n\}_{n\in\NN}$ such that $\cC_{k+1}(S'_1,\dots,S'_{k+1})\cap S'_{k+2}=\emptyset$ for any distinct $S'_1,\dots,S'_{k+2}$ chosen from $S_{i_1},\dots,S_{i_m}$, by Claim~\ref{claim4} we can find an $n_m\in\NN$ for which $Q(\vec{r},n_m)\cap \cC_{k+1}(S'_1,\dots,S'_{k+1})=\emptyset$.
    We pick an $S_{i_{m+1}}$ from $\{S_n\}_{n\in\NN}$ with $i_{m+1}>i_j$ for all $j\in [m]$ such that $S_{i_{m+1}}\subset Q(\vec{r},n_m)$.
    Continuing this way, we obtain our required $k$-independent sequence.
\end{proof}
Now it remains to prove the following claim:

\begin{claim}
    If every $k$-flat passing through $S_{i_1}, S_{i_2},\dots, S_{i_{k+1}}$ as defined in Lemma \ref{lemQ} does not contain any straight line parallel to $\vec{r}$, then $\exists n'\in\NN$ such that $Q(\vec{r},n')\cap\cC_{k+1}(S_{i_1}, S_{i_2},\dots, S_{i_{k+1}})=\emptyset$. 
\label{claim4}\end{claim}

\begin{proof}
    Let, without loss of generality, $\vec{r}=\{a e_d\mid a\geq 0\}$, where $e_d=(0,\dots,0,1)$.
    Since no $k$-flat through $S_{i_1}, S_{i_2},\dots, S_{i_{k+1}}$ contains any straight line parallel to $\vec{r}$, we conclude that through any $c_1,\dots,c_{k+1}$ with $c_j\in S_{i_j}$ for $j\in[k+1]$, exactly one $k$-flat passes.
    
    Since the sets $S_{i_1}, S_{i_2},\dots, S_{i_{k+1}}$ are compact and there is no $k$-flat parallel to $\vec{r}$ that intersects $S_{i_1}, S_{i_2},\dots, S_{i_{k+1}}$, we have a $c>0$ for which $\cC_{k+1}(S_{i_1}, S_{i_2},\dots, S_{i_{k+1}})\cap\{ae_d\in\RR^d\;|\;a\geq c \}=\emptyset$.
    This we can show in the following way: define a function $\tau: S_{i_1}\times\dots\times S_{i_{k+1}}\to\RR^{+}$ where $\tau(c_1,\dots,c_{k+1})$ is the distance of the $k$-flat that passes through $c_1,\dots,c_{k+1}$ from $\vec{r}$. 
    As $\tau$ is continuous, $\tau^{-1}(\{0\})$ is closed and therefore compact in $S_{i_1}\times\dots\times S_{i_{k+1}}$.
    Let $C=\tau^{-1}(\{0\})$.
    Then if $\rho:C\to\RR^{+}\cup\{0\}$ be the function for which $\rho(x_1,\dots,x_{k+1})=a$ implies that the $k$-flat that passes through $x_1,\dots,x_{k+1}$ intersects $\vec{r}$ at the point $ae_d$, then $\rho$ is continuous and therefore $\rho(C)$ is compact.
    From this, we get the required $c>0$.

    Now we define $\vec{r}'=ce_d+\vec{r}$ and $f: S_{i_1}\times S_{i_2}\times\dots\times S_{i_{k+1}}\rightarrow\RR^{+}$ by the following: take any $c_j\in S_j,\; \forall j\in [k+1]$ and suppose $\cK$ be the $k$-flat passing through $c_j\;\forall j\in [k+1]$ and $u_{\cK}\in\vec{r}'$ be the point closest to $\cK$. Then $$f(c_{i_1},\dots, c_{i_{k+1}})=\frac{\dist(\cK,\vec{r}')}{||u_{\cK}||}.$$ 
    For any given $k$-flat $\cK\in\cC_{k+1}(S_{i_1}, S_{i_2},\dots, S_{i_{k+1}})$, $u_{\cK}$ is uniquely defined as $\cK$ contains no straight line parallel to $\vec{r}$. 
    Also, $u_{\cK}\neq O$, as $O\not\in\vec{r}'$. 
    So, $f$ is well-defined. 
    Also, $\forall (c_{i_1},\dots, c_{i_{k+1}})\in S_{i_1}\times S_{i_2}\times\dots\times S_{i_{k+1}}$, $f(c_{i_1},\dots, c_{i_{k+1}})\neq 0$. 
    As $f$ is continuous, $\mathrm{image}(f)$ is compact. 
    Hence $\exists n'\in\NN$ such that $f(S_{i_1}\times S_{i_2}\times\dots\times S_{i_{k+1}})\cap [0,\frac{1}{n'}]=\emptyset$, and therefore, $\forall n>n'(1+\ceil{c})+1$, we have $Q(\vec{r},n)\cap\cC_{k+1}(S_{i_1}, S_{i_2},\dots, S_{i_{k+1}})=\emptyset$.
\end{proof}

\section{Heterochromatic variants of $(\aleph_0,q)$-theorems}\label{sec:heterochromatic}

\subsection{From monochromatic to heterochromatic}\label{ssec:mono_to_hetero}

Here we prove Theorem \ref{thm:mono_to_hetero}.
We begin with the following lemma:
\begin{lemma}\label{lem:hetero_bound}
    Let $\cF$ and $\cG$ be two collections of sets such that for every sequence $\{\cF_n\}_{n\in\NN}$ with $\cF_n\subseteq\cF$, $\{\cF_n\}_{n\in\NN}$ satisfying the heterochromatic $(\aleph_0,q)$-property with respect to $\cG$, there is an $n'\in\NN$ for which $\cF_{n'}$ has a finite-sized transversal $\cS_{n'}\subseteq\cG$.
    Then if $\{\cF_n\}_{n\in\NN}$ is a sequence with $\cF_n\subseteq\cF$ $\forall n\in\NN$ satisfying the heterochromatic $(\aleph_0,q)$-property with respect to $\cG$, then there is an $m\in\NN$ for which every $\cF_n$ with $n\geq m$ has a finite-sized transversal in $\cG$.
\end{lemma}

\begin{proof}
    If possible, let $\{\cF_{n_i}\}_{i\in\NN}$, where $\{n_i\}_{i\in\NN}$ is a strictly increasing sequence in $\NN$, be such that no $\cF_{n_i}$ has a finite-sized transversal in $\cG$.
    As $\{\cF_{n_i}\}_{i\in\NN}$ also satisfies the heterochromatic $(\aleph_0,q)$-property with respect to $\mathcal{G}$, there is an $i\in\NN$ for which $\cF_{n_i}$ has a finite-sized transversal in $\cG$.
    But this is a contradiction, and therefore, the result follows.
\end{proof}

Now we show that a finite transversal suffices for all but finitely many $\cF_n$s.

\begin{lemma}\label{lem:finite_fits_all}
    Let $\cF$ and $\cG$ be two collections of sets such that for every sequence $\{\cF_n\}_{n\in\NN}$ with $\cF_n\subseteq\cF$, $\{\cF_n\}_{n\in\NN}$ satisfying the heterochromatic $(\aleph_0,q)$-property with respect to $\cG$, there is an $n'\in\NN$ for which $\cF_{n'}$ has a finite-sized transversal $\cS_{n'}\subseteq\cG$.
    Then if $\{\cF_n\}_{n\in\NN}$ is a sequence with $\cF_n\subseteq\cF$ $\forall n\in\NN$ satisfying the heterochromatic $(\aleph_0,q)$-property with respect to $\cG$, then there is an $m\in\NN$ and a finite subset $\cS$ of $\cG$ such that $\cS$ is a transversal for all $\cF_n$ with $n\geq m$.
\end{lemma}

The difference between Lemma \ref{lem:finite_fits_all} and Lemma \ref{lem:hetero_bound} is that in Lemma \ref{lem:finite_fits_all}, we assert that there is one finite set that will work as a transversal for all but finitely many $\cF_n$s, whereas in Lemma \ref{lem:hetero_bound}, we show that all but finitely many $\cF_n$s have a finite-sized transversal in $\cG$.

\begin{proof}
    By Lemma \ref{lem:hetero_bound}, we have an $m\in\NN$ for which every $\cF_n$ with $n\geq m$ has a finite-sized transversal in $\cG$.
    Let $\cL=\bigcup_{n=m}^{\infty}\cF_n$.
    Then for any infinite subfamily $F\subseteq\cL$, either $F$ contains infinitely many elements of some $\cF_n$ with $n\geq m$, or $F$ contains a heterochromatic sequence with respect to the sequence $\{\cF_n\}_{n=m}^{\infty}$.
    In either case, this means that $F$ contains $q$ sets that are pierced by a set in $\cG$. 
    Therefore, $\cL$ satisfies the $(\aleph_0,q)$-property, and hence has a finite-sized transversal $\cS$ in $\cG$.
\end{proof}

Now, all that we need to prove Theorem \ref{thm:mono_to_hetero} is to show that a monochromatic $(\aleph_0,k+2)$-theorem for compact sets and $k$-flats implies its heterochromatic version, that is, in the terminology of Theorem \ref{thm:mono_to_hetero}, if every subset $\cF'$ of $\cF$ satisfying the $(\aleph_0,k+2)$-property with respect to $k$-flats has a finite $k$-transversal, then for every sequence $\{\cF_n\}_{n\in\NN}$ satisfying the heterochromatic $(\aleph_0,k+2)$-property with respect to $k$-flats, there is an $n'\in\NN$ for which $\cF_{n'}$ has a finite $k$-transversal.
For that, we shall need the following lemma, proved by Keller and Perles in \cite{KellerP22}.

\begin{lemma}\label{lem:KP}
    Let $\cF$ be a family of compact sets in $\RR^d$, $0\leq k\leq d-1$ and $m\in\NN$.
    If every finite subset of $\cF$ can be pierced by $m$ $k$-flats, then $\cF$ can be pierced by $m$ $k$-flats.
\end{lemma}

The next lemma will conclude the proof of Theorem \ref{thm:mono_to_hetero}.
\begin{lemma}\label{lem:mono_to_hetero}
    Let $d\in\NN$, $k\in\{0,\dots,d-1\}$ and $\cF$ be a collection of compact sets in $\RR^d$ such that if $\cF'\subseteq \cF$ satisfies the $(\aleph_0,k+2)$-property with respect to $k$-flats, then $\cF'$ has a finite $k$-transversal.
    Then if $\{\cF_n\}_{n\in\NN}$ is a sequence with $\cF_n\subseteq\cF$ satisfying the heterochromatic $(\aleph_0,k+2)$-property, then there is an $n'\in\NN$ for which $\cF_{n'}$ has a finite $k$-transversal.
\end{lemma}

\begin{proof}
    If for every $n\in\NN$, $\cF_n$ is not pierceable by finitely many $k$-flats, then by Lemma \ref{lem:KP}, we can obtain a finite subset $\cS_n\subset\cF_n$ such that $\cS_n$ is not pierceable by $n$ $k$-flats.
    Then $\cS:=\cup_{n\in\NN} \cS_n$ is not pierceable by finitely many $k$-flats, and therefore, there is an infinite subset $F\subseteq\cS$ that is $k$-independent.
    But as $F$ contains a heterochromatic sequence from $\{\cF_n\}_{n\in\NN}$, this violates the heterochromatic $(\aleph_0,k+2)$-property, and hence we arrive at a contradiction.
\end{proof}

\subsection{A strict heterochromatic $(\aleph_0,2)$-theorem}\label{ssec:strict_hetero}

Section \ref{sec:intro} and Subsection \ref{ssec:limits_strict_hetero} discuss how it is not possible to obtain a 'colorful' $(\aleph_0,k+2)$-theorem for compact convex sets in $\RR^d$ and $k$-flats with $k\in\{0,\dots,d-1\}$ if we only have the strict heterochromatic $(\aleph_0,k+2)$-property with respect to $k$-flats.
This is true even if we choose the sets to be closed unit balls in $\RR^d$ with $d>1$ and $k>0$.
However, we can prove Theorem \ref{thm:strong_hetero_balls}, which we do here.
\begin{proof}
	Pick $B_1\in\cB_1$.
	Try to pick a $B_2\in\cB_2$ that does not intersect $B_1$.
	If we can find one, then we move on to $\cB_3$ and try to pick a $B_3$ that does not intersect $B_1$ and $B_2$.
	Continuing this way, eventually we would find an $n_1\in\NN\setminus\{1\}$ for which there is no $B_{n_1}\in\cB_{n_1}$ that does not intersect $\cup_{i=1}^{n_1 -1} B_i$.
	Therefore, if $C=\{x\in\RR^d\mid \exists y\in\cup_{i=1}^{n_1 -1} B_i \text{such that } ||x-y||\leq 1\}$, then $\forall B_{n_1}\in\cB_{n_1}$, $B_{n_1}\subseteq C$.
	As $C$ is a compact set, it is easy to see that we can pierce all sets in $\cB_{n_1}$ by finitely many points.\\
	Let $B_0=\cup_{B\in\cB_{n_1}} B$.
	Now, let us try to pick a $B'_1\in\cB_1$ that does not intersect $B_0$.
	In general, for $i\in\NN$, let us try to pick a $B'_i\in\cB_i$ that does not intersect $\cup_{j=0}^{i-1}B_j$.
	There is an $n_2\in\NN\setminus\{n_1\}$ for which every $B\in\cB_{n_2}$ intersects $\cup_{j=0}^{n_2 -1}B_j$.
	Following the same arguments as in the case of $\cB_1$, we get that $\cB_{n_2}$ is pierceable by finitely many points.
	
\end{proof}

In the above theorem, we cannot have more than two distinct numbers, as is seen from the example where $\cB_{n_1}=\{B\}=\cB_{n_2}$ for some closed unit ball $B$, and every $\cB_i$ with $i\in\NN\setminus\{n_1,n_2\}$ is a $0$-set of closed unit balls.

\section{No-go results}\label{sec:no_go}

\subsection{Convexity and $\rho$-fatness in Theorem \ref{thm:mono_rho_fat}}

We begin by showing that the convexity assumption of Theorem \ref{thm:mono_rho_fat} cannot be removed without placing additional assumptions on the sets.

\begin{proposition}
	Given $\rho>0$, there is an infinite sequence $\{B_n\}_{n\in\NN}$ of compact connected $\rho$-fat sets in $\RR^2$ satisfying the $(\aleph_0,2)$-property with respect to points which cannot be pierced by finitely many points.
\end{proposition}

\begin{figure}
    \centering
    \includegraphics[height=5cm, width=15cm]{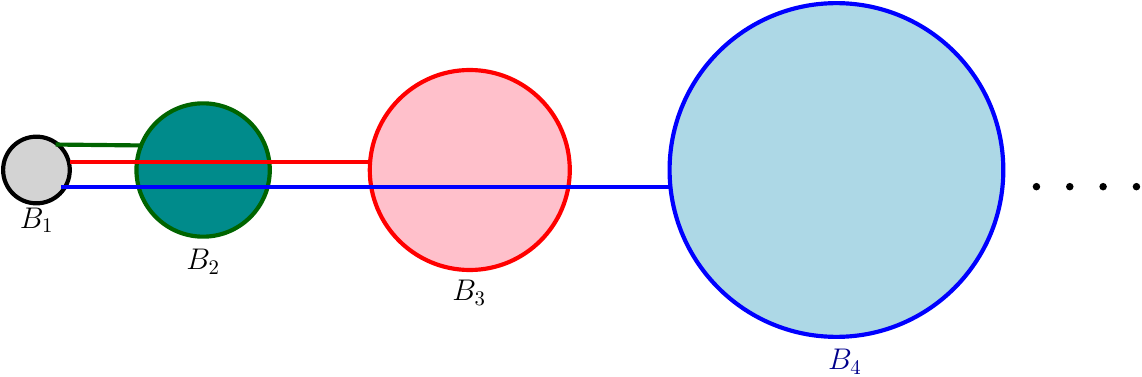}
    \caption{Sequence of non-convex sets with bounded condition number}
    \label{fig:impossibility_convexity}
\end{figure}

\begin{proof}
The construction here is similar to the construction in~\cite{KellerP22}.
First we pick a ball $B_1=B(c_1,r_1)$. Next we pick another ball $F_2=B(c_2,r_2)$ such that $B_1\subset B(c_2,\rho r_2)$ and $B_1\cap F_2=\emptyset$. Let $\Bar{\ell}_{2}$ be a compact curve joining $c_2$ and any point of $B_1$. Define $B_2=F_2\cup \Bar{\ell}_{2}$. In general, after $B_1, B_2,\dots, B_m$ are chosen, we pick another ball $F_{m+1}=B(c_{m+1},r_{m+1})$ such that $\cup_{i\in [m]}B_i\subset B(c_{m+1},\rho r_{m+1})$ and $\cup_{i\in [m]}B_i\cap F_{m+1}=\emptyset$.
Let $\Bar{\ell}_{m+1}$ be a compact curve from $c_{m+1}$ that passes through $B_i$ for every $i\in[m]$ and does not intersect $B_i\cap B_j$ for any $i,j\in[m]$, $i\neq j$.
Define $B_{m+1}=F_{m+1}\cup\Bar{\ell}_{m+1}$ (see Figure~\ref{fig:impossibility_convexity}). Then clearly, $\forall i\in [m],\; B_{m+1}\cap B_i\neq\emptyset$, i.e, $\{B_n\}_{n\in\NN}$ satisfies the $(\aleph_0,2)$ property. Again $\forall i,j\in [m]\text{ with } i\neq j,\; B_{m+1}\cap B_i\cap B_j=\emptyset$, i.e, no $3$ members of $\{B_n\}_{n\in\NN}$ have a common point. So $\{B_n\}_{n\in\NN}$ is not pierceable by a finite number of points.
\end{proof}

That the $\rho$-fat assumption too cannot be dropped without additional requirements on the sets can be seen from the fact that we can find an infinite collection of rectangles on the plane such that every pair of them intersect but the collection cannot be pierced by finitely many points.

\subsection{Limits of the strict heterochromatic $(\aleph_0,k+2)$-property}\label{ssec:limits_strict_hetero}

Proposition \ref{exmp:balls_hetero_no_go} shows that the strict heterochromatic $(\aleph_0,k+2)$-property is not enough to get a heterochromatic $(\aleph_0,k+2)$-theorem.
We now prove Proposition \ref{exmp:balls_hetero_no_go}.

\begin{proof}
	Pick a $k$-independent family $\cB_1=\{B_n\}_{n\in\NN}$ of closed unit balls in $\RR^d$.
	Let $\cB_i=\cB_1$ for $i=2,\dots,k+1$.
	Let $\cP = \{\alpha_1,\alpha_2,\dots\}$ denote the collection of all subsets of $\NN$ of size $k+1$.
	Let $S_{\alpha_j}$ denote the union of all $k$-flats that pass through each of $B_i$, $i\in\alpha$, $\alpha_j\in\cP$.
	In other words, $S_{\alpha_j}=\{x\in\RR^d\mid \text{there is a }k\text{-flat through }x\text{ that pierces each }B_i,i\in\alpha_j\}$.
	Now, 'tightly pack' $S_{\alpha_j}$ by closed unit balls.
	Denote the family of these balls by $\cB_{k+1+j}$.
	It is evident that $\cB_{k+1+j}$ is a $k$-set, since $S_{\alpha_j}$ contains an infinite cone 'packed' with closed unit balls, so for any finite collection of $k$-flats there is a closed unit ball in the cone that is far from those $k$-flats.

	Now, every strict heterochromatic sequence chosen from $\{\cB_n\}_{n\in\NN}$ contains $k+2$ sets that are pierced by a $k$-flat, because for any choice of sets $B_1,\dots,B_{k+1}$ from $\cB_1,\dots,\cB_{k+1}$ respectively, there is a $j\in\{k+2,k+3,\dots\}$ such that for any $B_j\in\cB_j$, there is a $k$-flat through $B_1,\dots,B_{k+1},B_j$.
	So, $\{\cB_n\}_{n\in\NN}$ satisfies the strict heterochromatic property with respect to $k$-flats.
	But, $\cB_n$ is a $k$-set for all $n\in\NN$.
\end{proof}

\subsection{An $(\aleph_0,q)$-theorem does not imply the existence of a $(p,q)$-theorem}\label{ssec:pq_no_go}

Here we prove Proposition \ref{exmp:near_balls_pq}.
\begin{proof}
	Choose $n$ closed balls $B_1,\dots,B_n$ in $\R^d$ with $B_i = B(c_i,r_i)$ such that for all $j\in\{2,\dots,n\}$, $\cup_{i=1}^{j-1}B(c_i,2r_i) \subset B(c_j,2r_j)$ and $\cup_{i=1}^{j-1}B(c_i,2r_i) \cap B_j =\emptyset$.
	Pick a set $P$ of $n\choose q$ distinct points in $B(c_1,2r_1)\setminus B_1$, where we write $P = \{p_S\mid S\text{ is a subset of size }q\text{ of }[n]\}$.
	Define, for each $i\in[n]$, $C_i = B_i\cup\{p_S\mid p_S\in P, i\in S\}$.
	Then $\{C_1,\dots,C_n\}$ is a family of compact near-balls, and for any $q$ distinct $i_1,\dots,i_q\in[n]$, there is a point $p_{\{i_1,\dots,i_q\}}\in P$ which belongs to every $C_{i_j}$ for $j\in[q]$.
	But, since every point in $\R^d$ intersects at most $q$ sets of this family, the smallest transversal for this family has size $\frac{n}{q}$.
	We can make every $C_i$ a contractible set by joining the points outside $B_i$ to $B_i$ with appropriate curves while maintaining the property that no point on $\R^d$ pierces more than $q$ sets of the new family.
\end{proof}

\section{Conclusion and open questions}

The connection between $(\aleph_0,q)$-theorems and their $(p,q)$-versions is relatively unexplored.
We know that a $(p,q)$-theorem and its corresponding fractional Helly theorem imply the $(\aleph_0,q)$-version.
We have also seen that we can get an $(\aleph_0,q)$-theorem even in the absence of its corresponding $(p,q)$-theorem.
However, it is still unknown under what conditions we can find a way to go from an $(\aleph_0,q)$-theorem to a $(p,q)$-theorem.
Future research throwing some light on this question would be interesting.

\bibliographystyle{alpha}
\bibliography{references}

\begin{thebibliography}{AKMM02}

\bibitem[AK92]{AlonK92}
N.~Alon and D.~J. Kleitman.
\newblock {Piercing Convex Sets and the Hadwiger-Debrunner $(p,q)$-Problem}.
\newblock {\em Advances in Mathematics}, 96(1):103--112, 1992.

\bibitem[AK95]{AlonK95}
N.~Alon and G.~Kalai.
\newblock {Bounding the Piercing Number}.
\newblock {\em Discrete \& Computational Geometry}, 13:245--256, 1995.

\bibitem[AKMM02]{alon2002transversal}
Noga Alon, Gil Kalai, Jiří Matou{\v{s}}ek, and Roy Meshulam.
\newblock Transversal numbers for hypergraphs arising in geometry.
\newblock {\em Advances in Applied Mathematics}, 29(1):79--101, 2002.

\bibitem[CGN25]{ChakrabortyGN2023boxesflats}
S.~Chakraborty, A.~Ghosh, and S.~Nandi.
\newblock {Stabbing boxes with finitely many axis-parallel lines and flats}.
\newblock {\em Discrete Mathematics}, 348(2):114269, 2025.

\bibitem[JP24a]{jung2024kdimensionaltransversalsfatconvex}
Attila Jung and Dömötör Pálvölgyi.
\newblock $k$-dimensional transversals for fat convex sets, 2024.

\bibitem[JP24b]{jung2024noteinfiniteversionspqtheorems}
Attila Jung and Dömötör Pálvölgyi.
\newblock A note on infinite versions of $(p,q)$-theorems, 2024.

\bibitem[KP22]{KellerP22}
C.~Keller and M.~A. Perles.
\newblock {An $(\aleph_{0}, k+2)$-Theorem for $k$-Transversals}.
\newblock In {\em Proceedings of the 38th International Symposium on
  Computational Geometry, {SoCG}}, volume 224, pages 50:1--50:14, 2022.

\bibitem[KP23]{KellerP23}
C.~Keller and M.~A. Perles.
\newblock {An $(\aleph_{0}, k+2)$-Theorem for k-Transversals}.
\newblock {\em CoRR}, abs/2306.02181, 2023.
\newblock To appear in {\em Israel Journal of Mathematics}.

\end{thebibliography}

\end{document}